\documentclass[11pt]{article}
\usepackage[top=1in,left=1in,right = 1in, bottom=0.7in, footskip = 0.2in]{geometry}
\usepackage[english]{babel}
\usepackage{algorithm,algorithmicx}
\usepackage{tcolorbox}
\usepackage{tikz}

\usepackage{framed}
\usepackage{algpseudocode}
\usepackage{amsthm,nccmath}
\usepackage{amsmath,bm}
\usepackage{bigints}
\usepackage{amssymb}
\usepackage{graphicx}
\usepackage{cite}
\usepackage[font=small,labelfont=bf]{caption}
\usepackage[font=scriptsize]{subcaption}
\usepackage{float}
\usepackage{epstopdf}
\usepackage{footnote}
\makesavenoteenv{tabular}
\makesavenoteenv{table}
\usepackage{bbm}
\usepackage{amsthm}
\usepackage[makeroom]{cancel}
\usepackage{stmaryrd}
\usepackage{amsmath,bm}
\usepackage{amssymb}
\usepackage{mathtools, cuted}
\usepackage{kantlipsum,setspace}

\usepackage{dblfloatfix}
\theoremstyle{plain}
\newtheorem{theorem}{Theorem}[section]
\newtheorem{lemma}[theorem]{Lemma}
\newtheorem*{lemma*}{Lemma}

\newtheorem*{cor*}{Corollary}

\theoremstyle{definition}

\newtheorem*{defn*}{Definition}

\theoremstyle{remark}
\newtheorem{rem}{Remark}
\newtheorem*{rem*}{Remark}

\newcommand{\nn}{\nonumber}
\newcommand{\beq}{\begin{equation}}
\newcommand{\eeq}{\end{equation}}
\newcommand{\bal}{\begin{align}}
\newcommand{\eal}{\end{align}}
\newcommand{\ba}{\begin{align*}}
\newcommand{\ea}{\end{align*}}
\newcommand{\bi}{\begin{itemize}}
\newcommand{\ei}{\end{itemize}}

\newcommand{\mbf}{\mathbf}

\newcommand{\lp}{\left(}
\newcommand{\rp}{\right)}

\newcommand{\lb}{\left[}
\newcommand{\rb}{\right]}
\newcommand{\lnr}{\left\|}
\newcommand{\rnr}{\right\|}

\usepackage{hyperref}

\usepackage{xcolor}

\DeclareMathOperator*{\argmin}{arg\,min}

\title{On Distributed Online Convex Optimization with Sublinear Dynamic Regret and Fit}

\author{Pranay Sharma$^{1}$, Prashant Khanduri$^{2}$, Lixin Shen$^{3}$, \\
Donald J. Bucci Jr.$^{4}$ and Pramod K. Varshney$^{1}$\\
{\small$^{1}$EECS, Syracuse University, Syracuse, NY}\\
{\small$^{2}$ECE, University of Minnesota, Minneapolis, MN}\\
{\small$^{3}$Mathematics, Syracuse University, Syracuse, NY}\\
{\small$^{4}$Lockheed Martin Advanced Technology Labs, Cherry Hill, NJ}\\
{\small Email: \{psharm04, lshen03, varshney\}@syr.edu,} \\
{\small khand095@umn.edu, donald.j.bucci.jr@lmco.com}
\thanks{P. Sharma, and P. K. Varshney are partially supported by the National Science Foundation (NSF) under grant ENG-1609916. This work was done when P. Khanduri was at Syracuse University. The research of L. Shen is partially supported by NSF under grant DMS-1913039.}
}

\begin{document}
\maketitle
\begin{abstract}
In this work, we consider a distributed online convex optimization problem, with time-varying (potentially adversarial) constraints. 
A set of nodes, jointly aim to minimize a global objective function, which is the sum of local convex functions. 
The objective and constraint functions are revealed locally to the nodes, at each time, after taking an action.
Naturally, the constraints cannot be instantaneously satisfied. 
Therefore, we reformulate the problem to satisfy these constraints in the long term. 
To this end, we propose a distributed primal-dual mirror descent-based approach.
% \PS{Should we have a line about what the algorithm does?}
The primal and dual updates are carried out locally at all the nodes. 
This is followed by sharing and mixing of the primal variables by the local nodes via communication with the immediate neighbors.
To quantify the performance of the proposed algorithm, we utilize the challenging, but more realistic metrics of dynamic regret and fit. 
Dynamic regret measures the cumulative loss incurred by the algorithm compared to the best dynamic strategy, while fit measures the long term cumulative constraint violations. 
Without assuming the restrictive Slater's conditions, we show that the proposed algorithm achieves sublinear regret and fit under mild, commonly used assumptions.
\end{abstract} 
% \begin{IEEEkeywords}
% Online optimization, Distributed optimization, Time-varying constraints, Dynamic regret, Mirror descent
% \end{IEEEkeywords}
% \IEEEpeerreviewmaketitle
 
\section{Introduction}
Many problems of practical interest, including network resource allocation \cite{chen17dynamic_tsp}, target tracking \cite{shahrampour17tae}, network routing \cite{zinkevich03icml}, online regression \cite{shalev12book}, and spam filtering \cite{hazan16book} can be framed in an Online Convex Optimization (OCO) framework. 
For a detailed review of OCO, please see \cite{shalev12book, hazan16book}. In this work, we consider
a constrained OCO problem, with time-varying (potentially
adversarial) constraints.
First introduced in \cite{zinkevich03icml}, the OCO framework aims to minimize a time varying convex objective function which is revealed to the observer in a sequential manner. 
% For a detailed review of OCO, please see \cite{hazan16book,shalev12book}. 
In this work, we consider a constrained OCO problem, with time-varying (potentially adversarial) constraints.

Recently, distributed OCO frameworks have gained popularity as they distribute the computations across multiple nodes rather than having a central node perform all the operations \cite{koppel15saddle_tsp, shahrampour17tae, ribeiro19static_dist_tv_constr, li18static_coupled, yi19dynamic_coupled}. 
% Distributed implementations make the algorithms robust by avoiding failure at a single central node, and also help enhance privacy.
We consider the constrained OCO problem in a distributed framework, where the convex objective is assumed to be decomposed and distributed across multiple communicating agents. 
Each agent takes its own action with the goal of minimizing the dynamically varying global function while satisfying its individual constraints.
%As pointed out later, the dynamic nature of the problem makes it infeasible to satisfy the constraints instantaneously, therefore, the best one can hope for is for the constraints to be satisfied in the long run.
Next, we discuss the related work along with the performance metrics we use to evaluate the performance of the proposed algorithm.

\subsection{Related Work}
\textbf{Regret:} The performance in OCO problems is quantified in terms of how well the agent does as compared to an offline system, over time. In other words, how much the agent ``regrets'' not having the information, which was revealed to it post-hoc, to begin with. Since regret is cumulative over time, an algorithm that achieves sub-linear increase in regret with time, asymptotically achieves zero average loss. It is naturally desirable to compare against an offline system, the action(s) of which are ``optimal'' in some sense.

\textbf{Static Regret:} The initial work on OCO, starting with \cite{zinkevich03icml, shalev12book, hazan16book}, almost exclusively focused on \textit{static regret} $\mathbf{Reg}_T^s$, which uses an optimal \textit{static} solution, in hindsight, as the benchmark. In other words, the fictitious offline adversary w.r.t. which the online system measures its regret, chooses the best fixed strategy, assuming it has access to the entire information, which is revealed to the online system over time horizon $T$.

\begin{align*}
    \mathbf{Reg}_T^s \triangleq \sum_{t=1}^T f_t (\mathbf{x}_t) - \min_{\mathbf{x}} \sum_{t=1}^T f_t (\mathbf{x}).
\end{align*}
Under standard regularity conditions, for general OCO problems, a tight upper bound of $O(\sqrt{T})$ has been shown for static regret \cite{zinkevich03icml, hazan07strong_convex}. However, for applications such as online parameter estimation or tracking moving targets, where the quantity of interest also evolves over time, comparison with a static benchmark is not sufficient.

This deficiency led to the development of \textit{dynamic regret} $\mathbf{Reg}_T^d$ \cite{hall15dynamic_jstsp, besbes15dynamic_or}. Rather than comparing the performance relative to a fixed optimal strategy, a more demanding benchmark is used. More precisely, at each time instant, our fictitious adversary utilizes one-step look-ahead information to adopt the optimal strategy at the current time instant.

\begin{align*}
    \mathbf{Reg}_T^d \triangleq \sum_{t=1}^T f_t (\mathbf{x}_t) - \sum_{t=1}^T \min_{\mathbf{x}} f_t (\mathbf{x}).
\end{align*}
In this work, we adopt the notion of dynamic regret as the performance metric. It must, however, be noted that, in the worst case, it is impossible to achieve sublinear dynamic regret \cite{zinkevich03icml}. For such problems, the growth of
dynamic regret is captured by the regularity measure which measures variations of the minimizer sequence over time (see $C_T^*$ in Theorem \ref{thm_regret_fit}).

\textbf{Constraints:} The conventional approaches for OCO are based on  projection-based gradient descent-like algorithms. 
However, when working with functional inequality constraints $g_t(\mathbf{x}) \leq \mathbf{0}$ (as opposed to simple convex feasible set constraints), the projection step in itself is computationally intensive.
This led to the development of primal-dual algorithms for OCO \cite{mahdavi12jmlr, jenatton16icml, yuan18cumul_nips}. 
Instead of attempting to satisfy the constraints at each time instant, the constraints are satisfied in the long run. 
In other words, the cumulative accumulation of instantaneous constraint violations (often simply called \textbf{fit}) $\| [\sum_{t=1}^T g_t (\mathbf{x}_t) ]_+ \|$ is shown to be sublinear in $T$. 
This formulation allows constraint violations at some instants to be ``taken-care-of'' by strictly feasible actions at other times.\footnote{Some more recent works \cite{yuan18cumul_nips} have considered the more stringent constraint violation metric $\sum_{t=1}^T ( [ g_t (\mathbf{x}_t) ]_+ )^2$}

Initially the constraints were assumed to be static across time \cite{mahdavi12jmlr, jenatton16icml}. 
However, subsequent literature \cite{sun17icml, chen17dynamic_tsp} demonstrated that the analysis for primal-dual methods can be generalized to even handle time-varying inequality constraints. 
Minor variations of primal-dual methods, which replace the dual update step with virtual-queue (modified Lagrange multiplier) updates have also been proposed to handle time-varying \cite{cao18dynamic_jstsp} and stochastic constraints \cite{neely17stoch_nips}.

\textbf{Distributed OCO Problems:} So far we have only discussed centralized problems. 
Suppose the OCO system has a network of agents, and local cost (and possibly constraint) functions are revealed to each agent over time. 
The global objective is to minimize the total cost function, while also satisfying all the constraints. And each agent can only communicate with those agents that are in its immediate neighborhood. 
This distributed OCO problem is more challenging and much less studied in the literature than the centralized problem.

Distributed OCO problems with static set constraints have been widely studied in recent years \cite{koppel15saddle_tsp, shahrampour17tae, ribeiro19static_dist_tv_constr, li18static_coupled, yi19dynamic_coupled}. 
Again here, the literature on distributed OCO with dynamic regret is much sparser than for static regret. 
The authors in \cite{shahrampour17tae} have proposed a dynamic mirror descent based algorithm, where primal update steps are alternated with local consensus steps. 
The authors in \cite{li18static_coupled} have proposed a distributed primal-dual algorithm for the OCO problem with coupled inequality constraints. 
The constraint functions are static over time. 
This has been generalized for time-varying coupled constraints in \cite{yi19dynamic_coupled}, where the authors have shown sublinearity of regret and fit, both w.r.t. dynamic and static benchmarks. 
However, to the best of our knowledge, the distributed OCO problem with a dynamic benchmark, even with static non-coupled inequality constraints has so far not been considered in the literature.

% Dynamic regret with constraints

% Distributed dynamic regret - with and without constraints

% Dynamic Regret papers - centralized without constraints - \cite{hall15dynamic_jstsp}; centralized with time-varying constraints - \cite{chen17dynamic_tsp}, \cite{cao18dynamic_jstsp}; distributed without constraints - \cite{shahrampour17tae};

% Coupled Inequality constraints - \cite{lee17static_coupled}, \cite{li18static_coupled}, \cite{yi19dynamic_coupled}

% Dist. time-varying constraints - static \cite{ribeiro19static_dist_tv_constr}

\subsection{Our Contributions} \label{subsec_contri}
In this work, we consider a distributed online convex optimization problem, where both the cost functions and the time-varying inequality constraints are revealed locally to the individual nodes.  We propose a primal-dual mirror-descent based algorithm, which alternates between the local primal and dual update steps and the consensus steps to mix the local primal variables with the immediate neighbors. Importantly, we show that the proposed algorithm achieves sublinear dynamic regret and fit.
%Our analysis has been borrowed in parts from \cite{shahrampour17tae} and \cite{yi19dynamic_coupled}.

\subsection{Paper Organization and Notations}
The paper is organized as follows: the problem formulation is discussed in Section \ref{sec_prob_form}, along with the definitions of the performance metrics. In Section \ref{sec_backgrnd}, we provide some background results and the assumptions required for providing theoretical guarantees. We propose our primal-dual mirror descent based algorithm in Section \ref{sec_alg}, followed by the theoretical results in Section \ref{sec_bounds}. Finally, we conclude the paper in Section \ref{sec_conc}.

\textbf{Notations:} Vectors are denoted with lowercase bold letters, e.g., $\mathbf{x}$, while matrices are denoted using uppercase bold letters, e.g., $\mathbf{X}$. The set of positive integers is represented by $\mathbb{N}_+$. We use $\mathbb{R}_+^n$ to denote the $n$-dimensional non-negative orthant. For $n \in \mathbb{N}_+$, the set $\{ 1,\hdots, n \}$ is denoted by $[n]$. We denote by $\| \cdot \|$ the Euclidean norm for vectors, and the induced 2-norm for matrices. $\mathbf{0}$ denotes a zero vector, where the dimension is clear from the context. $[ \mathbf{x} ]_+$ denotes the projection onto $\mathbb{R}_+^n$.

\section{Problem Formulation}
\label{sec_prob_form}
We consider a network of $n$ agents. At each time instant $t$, each agent $i$ takes an action $\mathbf{x}_{i,t} \in \mathcal{X} \subseteq \mathbb{R}^d$, where the set $\mathcal{X}$ is fixed across time, across all the nodes. Then, a set of local loss functions $\{ f_{i,t} (\cdot) \}_{i=1}^n$ with $f_{i,t} : \mathcal{X} \rightarrow \mathbb{R}$ are revealed to the individual nodes, leading to individual loss $f_{i,t} (\mathbf{x}_{i,t})$ at node $i$. Additionally, another set of local functions $\{ g_{i,t} (\cdot) \}_{i=1}^n$ with $g_{i,t}: \mathcal{X} \rightarrow \mathbb{R}^m$ are revealed, corresponding to local constraints $g_{i,t} (\mathbf{x}_{i,t}) \leq \mathbf{0}$. The network objective is to minimize the global average of the local cost functions $f_t (\mathbf{x}) \triangleq \tfrac{1}{n} \sum_{i=1}^n f_{i,t} (\mathbf{x})$, while also satisfying all the local constraint functions $\{ g_{i,t} (\cdot) \}_{i=1}^n$.

\beq
    \begin{aligned}
        & \min_{\mathbf{x}_t \in \mathcal{X}} \ f_t (\mathbf{x}_{t}) \triangleq \sum_{i=1}^n f_{i,t} (\mathbf{x}_{t}) \\
        & \text{subject to } g_{i,t} (\mathbf{x}_{t}) \leq \mathbf{0}_m, \forall \ i \in [n].
    \end{aligned}
    \label{eq_problem}
\eeq
Since the objective is to minimize the global function $f_t (\cdot)$, the nodes need to communicate among themselves. We next define the metrics used to measure the performance of the proposed approach.

\subsection{Performance Metrics - Dynamic Regret and Fit}
We use the notion of dynamic regret \cite{besbes15dynamic_or, hall15dynamic_jstsp} to measure the performance relative to a time-varying benchmark.

\begin{align}
    \mathbf{Reg}_T^d & \triangleq \frac{1}{n} \sum_{i=1}^n \sum_{t=1}^T f_t(\mathbf{x}_{i,t}) - \sum_{t=1}^T f_t(\mathbf{x}^*_t),  \label{eq_defn_dregret}
\end{align}
where $\mathbf{x}_{i,t}$ is the local action of agent $i$ at time $t$, while $\mathbf{x}^*_t$ is the solution of the following problem

\begin{align}
    \mathbf{x}^*_t \in \argmin_{\mathbf{x} \in \mathcal{X}} \left\{ f_t(\mathbf{x}) \big| g_{i,t} (\mathbf{x}) \leq \mathbf{0}, \forall \ i \in [n] \right\}.
\end{align}

As pointed out earlier, it is impossible to satisfy the time-varying constraints instantaneously, since they are revealed post-hoc. 
As a surrogate, to ensure the local constraints are satisfied in the \textit{long run}, we use the distributed extension of \textit{fit} as the performance metric. 
Fit has been used in the context of both time-invariant \cite{mahdavi12jmlr}, as well as time-varying constraints \cite{chen17dynamic_tsp, ribeiro16constr}, for single node problems. 
Our definition is motivated by the one given in \cite{ribeiro19static_dist_tv_constr} for continuous time problems. 
It measures the average accumulation of constraint violations over time.

\begin{align}
    \mathbf{Fit}_T^d & \triangleq \frac{1}{n} \sum_{i=1}^n \frac{1}{n} \sum_{j=1}^n \left\| \left[ \sum_{t=1}^T g_{i,t} (\mathbf{x}_{j,t}) \right]_+ \right\|. \label{eq_defn_dfit}
\end{align}
Here, $\sum_{t=1}^T g_{i,t} (\mathbf{x}_{j,t})$ is the constraint violation at agent $i$, if it adopts the actions of agent $j$. Note that $\sum_{t=1}^T g_{i,t} (\mathbf{x}_{j,t}) \leq \mathbf{0}$ is different from requiring the constraint to be met at every time instant $g_{i,t} (\mathbf{x}_{j,t}) \leq \mathbf{0}$.

Next, we discuss the assumptions and some background required for the analysis of the proposed OCO framework. Note that the following assumptions are standard for decentralized OCO problems \cite{shahrampour17tae, yi19dynamic_coupled}.

\section{Background and Assumptions}
\label{sec_backgrnd}
\subsection{Network}
\label{subsec_network}
We assume the $n$ agents are connected together via an undirected graph $\mathcal{G} = (\mathcal{V}, \mathcal{E})$. $\mathcal{V} = \{ 1,\hdots,n \}$ denotes the set of nodes of the graph, each of which represents an agent. $\mathcal{E}$ is the set of edges between the nodes. $(i,j) \in \mathcal{E}$ implies that nodes $i$ and $j$ are connected in the graph. The set of edges has an associated weight matrix $\mathbf{W}$, such that
\beq
    \mbf W = 
    \begin{cases}
        [\mathbf{W}]_{ij} > 0, & \qquad \text{if } (i,j) \in \mathcal{E}, \\
        [\mathbf{W}]_{ij} = 0 & \qquad \text{otherwise}.
    \end{cases}
\eeq
% $[\mathbf{W}]_{ij} > 0$ if  $(i,j) \in \mathcal{E}$, and  $[\mathbf{W}]_{ij} = 0$ otherwise.
The set of neighbors of node $i$ is, therefore, defined as $\mathcal{N}_i \triangleq \{ j : [\mathbf{W}]_{ij} > 0 \}$. Note that $j \in \mathcal{N}_i \Leftrightarrow i \in \mathcal{N}_j$.

\vspace{3mm}
\noindent \textbf{Assumption A:} \label{assum_nw} The network is connected. The weight matrix $\mathbf{W}$ is symmetric and doubly stochastic.
, such that
\begin{equation}
    \sum_{i=1}^n [\mathbf{W}]_{ij} = \sum_{j=1}^n [\mathbf{W}]_{ij} = 1.
\end{equation}
Next, we discuss the properties of the local cost functions and constraints.

\subsection{Local Objective Functions and Constraints}
\label{subsec_func}
\noindent \textbf{Assumption B:} \label{assum_func}
We assume the following conditions on the set $\mathcal{X}$, the objective and constraint functions.
\begin{enumerate}
    \item[(B1)] The set $\mathcal{X} \subseteq \mathbb{R}^d$ is convex and compact. Therefore, there exists a positive constant $d(\mathcal{X})$ such that
    \begin{equation}
        \left\| \mathbf{x} - \mathbf{y} \right\| \leq d(\mathcal{X}), \ \forall \ \mathbf{x}, \mathbf{y} \in \mathcal{\mathcal{X}}. \label{eq_set_diameter}
    \end{equation}
    \item[(B2)] The local node functions $f_{i,t}(\cdot), g_{i,t}(\cdot)$ are Lipschitz continuous on $\mathcal{X}$, $\forall \ i \in [n], \forall \ t \in \mathbb{N}_+$ i.e., for any $\mathbf{x}, \mathbf{y} \in \mathcal{X}$
    \begin{align}
        \begin{matrix}
            \left\| f_{i,t}(\mathbf{x}) - f_{i,t}(\mathbf{y}) \right\| \leq L \| \mathbf{x} - \mathbf{y} \|, \\
            \left\| g_{i,t}(\mathbf{x}) - g_{i,t}(\mathbf{y}) \right\| \leq L \| \mathbf{x} - \mathbf{y} \|.
        \end{matrix} \label{eq_defn_lipschitz}
    \end{align}
    \item[(B3)] The functions $\{ f_{i,t} \}, \{ g_{i,t} \}$ are convex and uniformly bounded on the set $\mathcal{X}$, i.e., there exists a constant $F > 0$ such that $\forall \ t \in \mathbb{N}_+, \forall \ i \in [n], \forall \ \mathbf{x} \in \mathcal{X}$
    \begin{equation}
        \| f_{i,t} (\mathbf{x}) \| \leq F, \| g_{i,t} (\mathbf{x}) \| \leq F.
    \end{equation}
    \item[(B4)] $\{ \nabla f_{i,t} \}, \{ \nabla g_{i,t} \}$ exist and are uniformly bounded on $\mathcal{X}$, i.e., there exists a constant $G > 0$ such that
    \begin{equation}
        \| \nabla f_{i,t} (\mathbf{x}) \| \leq G, \| \nabla g_{i,t} (\mathbf{x}) \| \leq G,
    \end{equation}
    $\forall \ t \in \mathbb{N}_+, \forall \ i \in [n], \forall \ \mathbf{x} \in \mathcal{X}$.
\end{enumerate}
Next, we briefly discuss Bregman Divergence, which is crucial to the proposed mirror descent based approach.

\subsection{Bregman Divergence}
\label{subsec_bregman}
Suppose we are given a $\mu$-strongly convex function $\mathcal{R}: \mathcal{X} \to \mathbb{R}$, i.e. $$\mathcal{R} (\mathbf{x}) \geq \mathcal{R} (\mathbf{y}) + \left\langle \nabla \mathcal{R} (\mathbf{y}), \mathbf{x} - \mathbf{y} \right\rangle + \frac{\mu}{2} \| \mathbf{x} - \mathbf{y} \|^2, \qquad \forall \ \mathbf{x}, \mathbf{y} \in \mathcal{X}.$$
The Bregman Divergence w.r.t. $\mathcal{R}$ is defined as
\begin{equation}
    \mathcal{D}_{\mathcal{R}} (\mathbf{x}, \mathbf{y}) \triangleq \mathcal{R} (\mathbf{x}) - \mathcal{R} (\mathbf{y}) - \left\langle \mathbf{x} - \mathbf{y}, \nabla \mathcal{R} (\mathbf{y}) \right\rangle. \label{eq_defn_bregman}
\end{equation}
Since $\mathcal{R} (\cdot)$ is $\mu$-strongly convex, for any $\mathbf{x}, \mathbf{y} \in \mathcal{X}$
\begin{equation}
    \mathcal{D}_{\mathcal{R}} (\mathbf{x}, \mathbf{y}) \geq \frac{\mu}{2} \| \mathbf{y} - \mathbf{x} \|^2. \label{eq_bregman_lower_bd}
\end{equation}
We assume the following conditions on $\mathcal{D}_{\mathcal{R}} (\cdot, \cdot)$.

\noindent \textbf{Assumption C:} \label{assum_bregman}
\begin{enumerate}
    \item[(C1)] Separate Convexity property \cite{bauschke01bregman}: Given $\mathbf{x}, \{\mathbf{y}_i \}_{i=1}^m \in \mathbb{R}^d$ and scalars $\{\alpha_i \}_{i=1}^m$ on the $m$-dimensional probability simplex, the Bregman Divergence satisfies
    
    \begin{align}
        \mathcal{D}_{\mathcal{R}} \left( \mathbf{x}, \sum_{i=1}^m \alpha_i \mathbf{y}_i \right) \leq \sum_{i=1}^m \alpha_i \mathcal{D}_{\mathcal{R}} \left( \mathbf{x}, \mathbf{y}_i \right). \label{eq_bregman_convex}
    \end{align}
    \item[(C2)] The Bregman divergence satisfies the following Lipschitz continuity condition \cite{jadbabaie15aistats}
    
    \begin{equation}
        \left| \mathcal{D}_{\mathcal{R}} \left( \mathbf{x}, \mathbf{y} \right) - \mathcal{D}_{\mathcal{R}} \left( \mathbf{z}, \mathbf{y} \right) \right| \leq K \left\|  \mathbf{x} - \mathbf{z} \right\| \label{eq_bregman_lipschitz}
    \end{equation}
    for any $\mathbf{x}, \mathbf{y}, \mathbf{z} \in \mathcal{X}$. This condition is satisfied if $\mathcal{R} (\cdot)$ is Lipschitz continuous on $\mathcal{X}$. Consequently,
    
    \begin{equation}
        \mathcal{D}_{\mathcal{R}} \left( \mathbf{x}, \mathbf{y} \right) \leq K d(\mathcal(X)), \ \forall \ \mathbf{x}, \mathbf{y} \in \mathcal{X}, \label{eq_bregman_bdd}
    \end{equation}
    where $d(\mathcal(X))$ is defined in \eqref{eq_set_diameter}.
\end{enumerate}

\noindent We next give a result on Bregman divergence from \cite{yi19dynamic_coupled} which is crucial to our analysis.

\begin{lemma}
\label{lemma_bregman}
Let $\mathcal{R}: \mathbb{R}^d \to \mathbb{R}$ be a $\mu$-strongly convex function. Also, assume $\mathcal{X}$ is a closed, convex set in $\mathbb{R}^p$ and $h: \mathcal{X} \to \mathcal{X}$ is a convex function. Assume that $\nabla h(\mathbf{x})$ exists $\forall \ \mathbf{x} \in \mathcal{X}$. Then, given $\mathbf{z} \in \mathcal{X}$, the regularized Bregman projection  $\mathbf{y} = \argmin_{\mathbf{x} \in \mathcal{X}} \left\{ h (\mathbf{x}) + \mathcal{D}_{\mathcal{R}} (\mathbf{x}, \mathbf{z}) \right\}$ satisfies the following inequality  $\forall \ \mathbf{x} \in \mathcal{X}$,

\begin{align}
    \left\langle \mathbf{y} - \mathbf{x}, \nabla h(\mathbf{y}) \right\rangle \leq \mathcal{D}_{\mathcal{R}} (\mathbf{x}, \mathbf{z}) - \mathcal{D}_{\mathcal{R}} (\mathbf{x}, \mathbf{y}) - \mathcal{D}_{\mathcal{R}} (\mathbf{y}, \mathbf{z}).
\end{align}
\end{lemma}

\subsection{Projection}
For a set $\mathcal{A} \subseteq \mathbb{R}^d$, the projection operator is defined as

\begin{align}
    \mathcal{P}_{\mathcal{A}} (\mathbf{y}) \triangleq \argmin_{\mathbf{x} \in \mathcal{A}} \left\| \mathbf{x} - \mathbf{y} \right\|^2, \label{eq_defn_proj}
\end{align}
$\forall \ \mathbf{y} \in \mathbb{R}^d$. For closed and convex $\mathcal{A}$, projection always exists and is unique. If $\mathcal{A} = \mathbb{R}_+^d$, projection is denoted by $[\cdot]_+$ and it satisfies

\begin{align}
    \left\| \left[ \mathbf{x} \right]_+ - \left[ \mathbf{y} \right]_+ \right\| \leq \left\| \mathbf{x} - \mathbf{y} \right\|, \forall \ \mathbf{x}, \mathbf{y} \in \mathbb{R}^d. \label{eq_proj_ineq}
\end{align}

\section{Distributed Primal-Dual Mirror Descent based Algorithm}
\label{sec_alg}
We next discuss the proposed distributed primal-dual mirror descent based algorithm for online convex optimization with time-varying constraints. The pseudo-code is outlined in Algorithm \ref{alg1}. The algorithm runs in parallel at all the nodes. At the end of time $t-1$, $\mathbf{x}_{i,t-1}$ is the action (primal variable) at node $i$. Following this, the local functions $f_{i,t-1}, g_{i,t-1}$ are revealed to the agent. The corresponding function values and gradients are utilized to carry-out the updates in the next time step $t$. First, each agent performs the primal update locally  (Step \ref{alg1_line_primal_update}). This is followed by the dual update (Step \ref{alg1_line_dual_update}). Note that the projection $[\cdot]_+$ ensures that the dual variable lies in the non-negative orthant $\mathbb{R}_+^m$. At the end of each time step, an average consensus step is taken across the nodes, where the local updated primal variables $\mathbf{y}_{i,t-1}$ are received from the neighbors, to compute the action $\mathbf{x}_{i,t}$.

\begin{algorithm}[t!]
\caption{Distributed OCO with constraints}
\label{alg1}
\begin{algorithmic}[1]
	\State{\textbf{Input}: Non-increasing sequences $\{ \alpha_t >0 \}, \{ \beta_t >0 \}, \{ \gamma_t >0 \}$; Differentiable and strongly-convex $\mathcal{R}$}
	\State{\textbf{Initialize}: $\mathbf{x}_{i,0} = \mathbf{0}_d \in \mathcal{X}$, $f_{i,0}(\cdot) \equiv 0, g_{i,0}(\cdot) \equiv \mathbf{0}_m$, $\mathbf{q}_{i,0} = \mathbf{0}_m$, $\forall \ i \in [n]$}.
	\For{$t = 1$ to $T$}
    	\For{$i = 1$ to $n$}
    	    \State{Observe $\nabla f_{i,t-1}(\mathbf{x}_{i,t-1}), \nabla g_{i,t-1}(\mathbf{x}_{i,t-1}), g_{i,t-1}(\mathbf{x}_{i,t-1})$}
        	\State{\label{alg_line_update_a} $\mathbf{a}_{i,t} = \nabla f_{i,t-1}(\mathbf{x}_{i,t-1}) + \left[ \nabla g_{i,t-1}(\mathbf{x}_{i,t-1}) \right]^T \mathbf{q}_{i,t-1}$}
        	\State{\label{alg1_line_primal_update} $\mathbf{y}_{i,t} = \argmin_{\mathbf{x} \in \mathcal{X}} \left\{ \alpha_t \left\langle \mathbf{x}, \mathbf{a}_{i,t} \right\rangle + \mathcal{D}_{\mathcal{R}} (\mathbf{x}, \mathbf{x}_{i,t-1}) \right\}$ \hspace*{\fill} (Local Primal update)}
        	\State{\label{alg_line_update_b} $\mathbf{b}_{i,t} = \left[ \nabla g_{i,t-1}(\mathbf{x}_{i,t-1}) \right] ( \mathbf{y}_{i,t} - \mathbf{x}_{i,t-1} ) + g_{i,t-1}(\mathbf{x}_{i,t-1})$}
        	\State{\label{alg1_line_dual_update} $\mathbf{q}_{i,t} = \left[ \mathbf{q}_{i,t-1} + \gamma_t (\mathbf{b}_{i,t} - \beta_{t} \mathbf{q}_{i,t-1}) \right]_+$ \hspace*{\fill} (Local dual update)}
        	\State{Broadcast $\mathbf{y}_{i,t}$ to out-neighbors $j \in \mathcal{N}_i$}
        	\State{Obtain $\mathbf{y}_{j,t}$ from in-neighbors $j \in \mathcal{N}_i$}
        	\State{\label{alg1_line_consensus} $\mathbf{x}_{i,t} = \sum_{j=1}^n [\mathbf{W}]_{ij} \mathbf{y}_{j,t}$ \hspace*{\fill} (Consensus step)}
	    \EndFor
	\EndFor
\end{algorithmic}
\end{algorithm}

\begin{rem}
Note that the primal and dual update steps employ different step-sizes, $\alpha_t$ and $\gamma_t$, respectively. This idea originated in \cite{jenatton16icml} and leads to flexibility in terms of the trade-off between the bounds on dynamic regret and fit.
\end{rem}
In the next section, we bound the dynamic regret and fit which result from Algorithm \ref{alg1}, and show them to be sublinear in the time-horizon $T$.

\section{Dynamic Regret and Fit Bounds}
\label{sec_bounds}
First, we discuss some intermediate results required to show the sublinearity of dynamic regret and fit. 
We have omitted the proofs due to space limitations. 
Our analysis follows closely the work in \cite{shahrampour17tae} and \cite{yi19dynamic_coupled}.
\subsection{Some Intermediate Results}
\label{subsec_interm_lemmas}
\begin{lemma}
\label{lemma_dual}
Suppose Assumption B holds. $\forall \ i \in [n]$, $\forall \ t \in \mathbb{N}_+$, $\mathbf{q}_{i,t}$ generated by Algorithm \ref{alg1} satisfy

\begin{align}
    \left\| \mathbf{q}_{i,t} \right\| & \leq \frac{F}{\beta_{t}} \label{eq_lemma_dual_1} \\
    \frac{\Delta_{t+1}}{2 \gamma_{t+1}} & \leq \frac{n B_1^2}{2} \gamma_{t+1} +  \sum_{i=1}^n \mathbf{q}_{i,t}^T [\nabla g_{i,t} (\mathbf{x}_{i,t})] (\mathbf{y}_{i,t+1} - \mathbf{x}_{i,t}) + \left( \frac{G^2 \alpha_{t+1}}{\mu} + \frac{\beta_{t+1}}{2} \right) \sum_{i=1}^n \| \mathbf{q}_{i} \|^2 \nn \\
    & \quad + \sum_{i=1}^n (\mathbf{q}_{i,t} - \mathbf{q}_{i})^T g_{i,t} (\mathbf{x}_{i,t}) + \frac{\mu}{4 \alpha_{t+1}} \sum_{i=1}^n \| \mathbf{y}_{i,t+1} - \mathbf{x}_{i,t} \|^2 \label{eq_lemma_dual_2}
\end{align}
where $B_1 = 2F + G d(\mathcal{X})$,

\begin{align*}
    \Delta_{t+1} \triangleq \sum_{i=1}^n \left[ \| \mathbf{q}_{i,t+1} - \mathbf{q}_{i} \|^2 - (1 - \gamma_{t+1} \beta_{t+1}) \| \mathbf{q}_{i,t} - \mathbf{q}_{i} \|^2 \right],
\end{align*}
and $\{ \mathbf{q}_i \}_i$ are arbitrary vectors in $\mathbb{R}^m_+$.
\end{lemma}

% \begin{proof}
% See Appendix \ref{app_lemma_dual}.
% \end{proof}

\begin{rem}
The penalty term $- \beta_{t} \mathbf{q}_{i,t-1}$ in the dual update (step \ref{alg1_line_dual_update}, Algorithm \ref{alg1}) helps in upper bounding the local dual variables. This idea was initially used in \cite{mahdavi12jmlr} and helps get rid of the requirement of Slater's condition. $\Delta_{t+1}$ measures the regularized drift of the local dual variables. See \cite{hall15dynamic_jstsp} and \cite{yi19dynamic_coupled} for similar results, respectively in centralized and distributed contexts.
\end{rem}

Next, we sum the left hand side of \eqref{eq_lemma_dual_2} over $t$ to get
\beq
    \begin{aligned}
        \sum_{t=1}^T \frac{\Delta_{t+1}}{2 \gamma_{t+1}} &= \frac{1}{2} \sum_{t=1}^T \left( \frac{1}{\gamma_t} - \frac{1}{\gamma_{t+1}} + \beta_{t+1} \right) \sum_{i=1}^n \| \mathbf{q}_{i,t} - \mathbf{q}_{i} \|^2 \\
        & \quad -\frac{1}{2} \sum_{i=1}^n \left[ \frac{1}{\gamma_1} \| \mathbf{q}_{i,1} - \mathbf{q}_{i} \|^2 - \frac{1}{\gamma_{T+1}} \| \mathbf{q}_{i,T+1} - \mathbf{q}_{i} \|^2 \right]. \label{eq_Delta_sum_over_t}
    \end{aligned}
\eeq
Recall that $\mathbf{q}_{i,1} = \mathbf{0}$, $\forall \ i \in [n]$. We combine \eqref{eq_lemma_dual_2} and \eqref{eq_Delta_sum_over_t}, and define $g_c(\cdot)$ such that

\begin{align}
    & g_c (\mathbf{q}_1, \hdots, \mathbf{q}_n) \triangleq \sum_{i=1}^n \mathbf{q}_{i}^T \Big( \sum_{t=1}^T g_{i,t} (\mathbf{x}_{i,t}) \Big) - \left[ \frac{1}{2 \gamma_1} + \sum_{t=1}^T \left( \frac{G^2 \alpha_{t+1}}{\mu} + \frac{\beta_{t+1}}{2} \right) \right] \sum_{i=1}^n \| \mathbf{q}_{i} \|^2 \nonumber \\
    & \quad \leq \frac{n B_1^2}{2} \sum_{t=1}^T \gamma_{t+1} + \sum_{t=1}^T \sum_{i=1}^n \mathbf{q}_{i,t}^T [\nabla g_{i,t} (\mathbf{x}_{i,t})] (\mathbf{y}_{i,t+1} - \mathbf{x}_{i,t}) + \sum_{t=1}^T \sum_{i=1}^n \mathbf{q}_{i,t}^T g_{i,t} (\mathbf{x}_{i,t}) \nn \\
    & \qquad + \sum_{t=1}^T \frac{\mu}{4 \alpha_{t+1}} \sum_{i=1}^n \| \mathbf{y}_{i,t+1} - \mathbf{x}_{i,t} \|^2 - \frac{1}{2} \sum_{t=1}^T \left( \frac{1}{\gamma_t} - \frac{1}{\gamma_{t+1}} + \beta_{t+1} \right) \sum_{i=1}^n \| \mathbf{q}_{i,t} - \mathbf{q}_{i} \|^2. \label{eq_bd_gc_q}
\end{align}
The function $g_c (\mathbf{q}_1, \ldots, \mathbf{q}_n)$ will be used later in Lemma \ref{lemma_regret_fit} to upper bound both the dynamic regret and fit, by appropriately choosing $\mathbf{q}_i, \ \forall \ i \in [n]$.
Before looking at the primal updates, we first consider one of the constituent terms in \eqref{eq_defn_dregret}.
\begin{align}
    f_t(\mathbf{x}_{i,t}) - f_t(\mathbf{x}^*_t) &= f_t(\mathbf{x}_{i,t}) - f_t(\bar{\mathbf{x}}_{t}) + f_t(\bar{\mathbf{x}}_{t}) - f_t(\mathbf{x}^*_t) \nonumber \\
    & \leq L \left\| \mathbf{x}_{i,t} - \bar{\mathbf{x}}_{t} \right\| + f_t(\bar{\mathbf{x}}_{t}) - f_t(\mathbf{x}^*_t) \label{eq_use_f_lipschitz_1} \\
    & = \frac{1}{n} \sum_{j=1}^n \left\{ f_{j,t} (\bar{\mathbf{x}}_{t}) - f_{j,t} (\mathbf{x}^*_t) + f_{j,t} (\mathbf{x}_{j,t}) - f_{j,t} (\mathbf{x}_{j,t}) \right\} + L \left\| \mathbf{x}_{i,t} - \bar{\mathbf{x}}_{t} \right\| \nonumber \\
    & \leq \frac{1}{n} \sum_{j=1}^n \left\{ f_{j,t} (\mathbf{x}_{j,t}) - f_{j,t} (\mathbf{x}^*_t) \right\} + L \left\| \mathbf{x}_{i,t} - \bar{\mathbf{x}}_{t} \right\| + \frac{L}{n} \sum_{j=1}^n \left\| \mathbf{x}_{j,t} - \bar{\mathbf{x}}_{t} \right\|. \label{eq_use_f_lipschitz_2}
\end{align}
We use assumption (B2) to obtain both \eqref{eq_use_f_lipschitz_1}, \eqref{eq_use_f_lipschitz_2}.
Now, from the definition of dynamic regret \eqref{eq_defn_dregret}, we get
\begin{align}
    \mathbf{Reg}_T^d & \leq \frac{1}{n} \sum_{i=1}^n \sum_{t=1}^T \frac{1}{n} \sum_{j=1}^n \left\{ f_{j,t} (\mathbf{x}_{j,t}) - f_{j,t} (\mathbf{x}^*_t) \right\} + \frac{2L}{n} \sum_{i=1}^n \sum_{t=1}^T \left\| \mathbf{x}_{i,t} - \bar{\mathbf{x}}_{t} \right\|. \label{eq_bd_dregret_11}
\end{align}
Next, we upper bound both the terms in \eqref{eq_bd_dregret_11}. First, we upper bound the first term in the following lemma.

\begin{lemma}
\label{lemma_primal}
Suppose Assumptions A-C hold. $\forall \ i \in [n]$, $\forall \ t \in \mathbb{N}_+$, if $\{ \mathbf{x}_{i,t} \}$ is the sequence generated by Algorithm \ref{alg1}. Then,
\beq
    \begin{aligned}
        & \sum_{t=1}^T \sum_{i=1}^n \left[ f_{i,t} (\mathbf{x}_{i,t}) - f_{i,t} (\mathbf{x}^*_{t}) \right] \leq \frac{n G^2}{\mu} \sum_{t=1}^T \alpha_{t+1} - \sum_{t=1}^T \sum_{i=1}^n \frac{\mu}{4 \alpha_{t+1}} \| \mathbf{y}_{i,t+1} - \mathbf{x}_{i,t} \|^2 \\
        & \qquad \qquad + \frac{n K}{\alpha_{T+2}} \sum_{t=1}^T \| \mathbf{x}^*_{t+1} - \mathbf{x}^*_{t} \| - \sum_{t=1}^T \sum_{i=1}^n \mathbf{q}_{i,t}^T \Big[ g_{i,t}(\mathbf{x}_{i,t}) + \nabla g_{i,t}(\mathbf{x}_{i,t}) \left( \mathbf{y}_{i,t+1} - \mathbf{x}_{i,t} \right) \Big] \\
        & \qquad \qquad + \sum_{i=1}^n \left[ \frac{1}{\alpha_{2}} \mathcal{D}_{\mathcal{R}} (\mathbf{x}^*_{1}, \mathbf{x}_{i,1}) - \frac{1}{\alpha_{T+2}} \mathcal{D}_{\mathcal{R}} \lp \mathbf{x}^*_{T+1}, \mathbf{x}_{i,T+1} \rp \right] + \frac{n K d(\mathcal(X))}{\alpha_{T+2}}.
    \end{aligned}
    \label{eq_lemma_dregret}
\eeq
\end{lemma}

% \begin{proof}
% See Appendix \ref{app_lemma_primal}.
% \end{proof}
Next, we upper bound the second term in \eqref{eq_bd_dregret_11}. This is the consensus error of the primal variables.
\begin{lemma}{(Network Error):}
\label{lemma_nw_error}
Suppose Assumptions A-C hold. Then, the local estimates $\{ \mathbf{x}_{i,t}\}$ generated by Algorithm \ref{alg1} satisfy
\begin{align}
    \left\| \mathbf{x}_{i,t} - \bar{\mathbf{x}}_t \right\| \leq \sum_{\tau = 0}^{t-1} \sqrt{n} \sigma_2^{t-\tau} (\mathbf{W}) \frac{G \alpha_{\tau+1}}{\mu} \left( 1 + \frac{F}{\beta_{\tau+1}} \right) \label{eq_nw_error}
\end{align}
$\forall \ i \in [n]$, where $\bar{\mathbf{x}}_t = \frac{1}{n} \sum_{i=1}^n \mathbf{x}_{i,t}$. $\sigma_2 (\mathbf{W})$ is the second largest eigenvalue of $\mathbf{W}$ in magnitude.
\end{lemma}

\begin{rem}
The network error bound is \eqref{eq_nw_error} is independent of the node index $i$. The dependence on $\sigma_2 (\mathbf{W})$ captures the speed with which mixing of iterates happens. The smaller the value of $\sigma_2 (\mathbf{W})$, the faster the network error diminishes. Moreover, the choice of the primal update step sizes $\{ \alpha_t \}$ and the dual update regularization parameters $\{ \beta_t \}$ has a crucial role to play in bounding the network error. As we shall see in Theorem \ref{thm_regret_fit}, carefully choosing these leads to sublinear regret and fit.
\end{rem}

% \begin{proof}
% See Appendix \ref{app_lemma_nw_error}.
% \end{proof}

Next, we combine \eqref{eq_bd_gc_q} and Lemma \ref{lemma_primal} resulting in two intermediate bounds, which shall be needed to subsequently bound the dynamic regret and fit respectively.
\begin{lemma}
\label{lemma_regret_fit}
Suppose Assumptions A-C hold. Then, the sequences $\{ \mathbf{x}_{i,t}, \mathbf{q}_{i,t} \}$ generated by Algorithm \ref{alg1} satisfy
\beq
    \begin{aligned}
        & \sum_{t=1}^T \sum_{i=1}^n \left( f_{i,t} (\mathbf{x}_{i,t}) - f_{i,t} (\mathbf{x}^*_{t}) \right) \leq \frac{n B_1^2}{2} \sum_{t=1}^T \gamma_{t+1} + \frac{n G^2}{\mu} \sum_{t=1}^T \alpha_{t+1} + \frac{n K}{\alpha_{T+2}} \sum_{t=1}^T \lnr \mathbf{x}^*_{t+1} - \mathbf{x}^*_{t} \rnr \\
        & \qquad \qquad  + \frac{n K d(\mathcal(X))}{\alpha_{T+2}} + \sum_{i=1}^n \left[ \frac{1}{\alpha_{2}} \mathcal{D}_{\mathcal{R}} (\mathbf{x}^*_{1}, \mathbf{x}_{i,1}) - \frac{1}{\alpha_{T+2}} \mathcal{D}_{\mathcal{R}} (\mathbf{x}^*_{T+1}, \mathbf{x}_{i,T+1}) \right] \\
        & \qquad \qquad - \frac{1}{2} \sum_{t=1}^T \left( \frac{1}{\gamma_t} - \frac{1}{\gamma_{t+1}} + \beta_{t+1} \right) \sum_{i=1}^n \| \mathbf{q}_{i,t} \|^2, \label{eq_lemma_dregret_2}
    \end{aligned}
\eeq
and
\beq
    \begin{aligned}
        & \sum_{i=1}^n \left\| \left[ \sum_{t=1}^T g_{i,t} (\mathbf{x}_{i,t}) \right]_+ \right\|^2 \leq 4\left[ \frac{1}{2 \gamma_1} + \sum_{t=1}^T \left( \frac{G^2 \alpha_{t+1}}{\mu} + \frac{\beta_{t+1}}{2} \right) \right] \Bigg\{ 2nFT \\
        & \qquad + \frac{n B_1^2}{2} \sum_{t=1}^T \gamma_{t+1} + \frac{n G^2}{\mu} \sum_{t=1}^T \alpha_{t+1}  + \sum_{i=1}^n \left[ \frac{1}{\alpha_{2}} \mathcal{D}_{\mathcal{R}} (\mathbf{x}^*_{1}, \mathbf{x}_{i,1}) - \frac{1}{\alpha_{T+2}} \mathcal{D}_{\mathcal{R}} (\mathbf{x}^*_{T+1}, \mathbf{x}_{i,T+1}) \right] \\
        & \qquad + \frac{n K}{\alpha_{T+2}} \sum_{t=1}^T \| \mathbf{x}^*_{t+1} - \mathbf{x}^*_{t} \| + \frac{n K d(\mathcal(X))}{\alpha_{T+2}}  - \frac{1}{2} \sum_{t=1}^T \left( \frac{1}{\gamma_t} - \frac{1}{\gamma_{t+1}} + \beta_{t+1} \right) \sum_{i=1}^n \| \mathbf{q}_{i,t} - \bar{\mathbf{q}}_i \|^2 \Bigg\}.
    \end{aligned}
    \label{eq_lemma_dfit}
\eeq
\end{lemma}

\begin{rem}
\eqref{eq_lemma_dregret_2} follows by adding \eqref{eq_bd_gc_q} and \eqref{eq_lemma_dregret}, and substituting $\mathbf{q}_i = \mathbf{0}_m, \ \forall \ i \in [n]$. Similarly, \eqref{eq_lemma_dfit} is obtained by adding \eqref{eq_bd_gc_q} and \eqref{eq_lemma_dregret}, and substituting

\begin{align}
    \bar{\mathbf{q}}_i = \frac{\left[ \sum_{t=1}^T g_{i,t} (\mathbf{x}_{i,t}) \right]_+}{ \gamma_1^{-1} + \sum_{t=1}^T \left( 2 G^2 \alpha_{t+1}/\mu + \beta_{t+1} \right) }, \ \forall \ i \in [n].
\end{align}
\end{rem}

% \begin{proof}
% See Appendix \ref{app_lemma_regret_fit}.
% \end{proof}
Before presenting out final result, we need to use the following upper bound to bound the fit.
\begin{align}
    \frac{1}{n} \sum_{i=1}^n \frac{1}{n} \sum_{j=1}^n \Big\| \left[ \sum_{t=1}^T g_{i,t} (\mathbf{x}_{j,t}) \right]_+ \Big\|^2 \leq 2 \left[ 2L \sum_{t=1}^T \left\| \mathbf{x}_{i,t} - \bar{\mathbf{x}}_t \right\| \right]^2 + \frac{2}{n} \sum_{i=1}^n \lnr \left[ \sum_{t=1}^T g_{i,t} (\mathbf{x}_{i,t}) \right]_+ \rnr^2. \label{eq_fit_sq}
\end{align}
This follows from Lipschitz continuity of the constraint functions (Assumption (B2)). Since, we have bounded both the terms in \eqref{eq_fit_sq} (the first term in Lemma \ref{lemma_nw_error}, and the second term in Lemma \ref{lemma_regret_fit}), we are now ready to present our final result on the sublinearity of both dynamic regret and fit.
 
\subsection{Dynamic Regret and Fit Bounds}
\begin{theorem}
\label{thm_regret_fit}
Suppose Assumptions A-C hold, and $\{ \mathbf{x}_{i,t}\}$ be the sequence of local estimates generated by Algorithm \ref{alg1}. We choose the step sizes
\begin{equation}
    \alpha_t = \frac{1}{t^{a}}, \ \beta_t = \frac{1}{t^{b}}, \  \gamma_t = \frac{1}{t^{1-b}}, \quad \forall \ t \in \mathbb{N}_+ \label{eq_step_size_choice}
\end{equation}
where, $a, b \in (0,1)$ and $a > b$. Then for any $T \in \mathbb{N}_+$.

\begin{align}
    \mathbf{Reg}_T^d & \leq R_{1} T^{\max \{ a, 1-a+b \}} + 2 K T^a C_T^*, \label{eq_thm1_regret} \\
    \frac{1}{n} \sum_{i=1}^n \frac{1}{n} \sum_{j=1}^n \left\| \left[ \sum_{t=1}^T g_{i,t} (\mathbf{x}_{j,t}) \right]_+ \right\|^2 & \leq D_1 T^{2-b} + D_2 T^{1+a-b} C_T^* + D_3 T^{2+2b-2a} \label{eq_thm1_fit}.
\end{align}
Here, $R = \frac{4 F L G \sqrt{n} \sigma_2(\mathbf{W})}{\mu (1-a) (1 - \sigma_2(\mathbf{W}))}, R_1 = R + \frac{B_1^2}{2b} + \frac{G^2}{\mu (1-a)} + 2 K d(\mathcal(X)), D = 2 + \frac{4 G^2}{\mu (1-a)} + \frac{2}{1-b} , D_1 = 2 D (2F + 2 K d(\mathcal(X)) + \frac{B_1^2}{2 b} + \frac{G^2}{\mu (1-a)} + 2 K d(\mathcal(X)))$, $D_2 = 4 K D$ and $D_3 = 16L^2 R^2$ are constants independent of $T$, and

\begin{equation}
    C_T^* \triangleq \sum_{t=1}^T \| \mathbf{x}^*_{t+1} - \mathbf{x}^*_{t} \|
\end{equation}
is the accumulated dynamic variation of the comparator sequence $\{ \mathbf{x}_t^* \}$.
\end{theorem}

% \begin{proof}
% See Appendix \ref{app_thm_regret_fit}.
% \end{proof}

\begin{rem}
The dynamic regret $\mathbf{Reg}_T^d$ is sublinear as long as the cumulative consecutive variations of the dynamic comparators $C_T^*$ is sublinear. This is the standard requirement for sublinearity of dynamic regret \cite{hall15dynamic_jstsp, shahrampour17tae, yi19dynamic_coupled}.
\end{rem}

\begin{rem}
A similar argument as above holds for \eqref{eq_thm1_fit}. As long as $C_T^*$ is sublinear, we have
\begin{align}
    \frac{1}{n} \sum_{i=1}^n \frac{1}{n} \sum_{j=1}^n \left\| \left[ \sum_{t=1}^T g_{i,t} (\mathbf{x}_{j,t}) \right]_+ \right\|^2 = o(T^2). \label{eq_fit_sq2}
\end{align}
Note that \eqref{eq_thm1_fit} has $\| [ \sum_{t=1}^T g_{i,t} (\mathbf{x}_{j,t}) ]_+ \|^2$, while fit \eqref{eq_defn_dfit} is defined with $\| [ \sum_{t=1}^T g_{i,t} (\mathbf{x}_{j,t}) ]_+ \|$. However, for large enough $T$, each of the constituent terms in \eqref{eq_fit_sq2} are $o(T^2)$. Consequently, 
$$\| [ \sum_{t=1}^T g_{i,t} (\mathbf{x}_{j,t}) ]_+ \|^2 = o(T),$$
$\forall \ i, j \in [n]$. Therefore, we get a sublinear fit

\begin{align}
    \mathbf{Fit}_T^d = \frac{1}{n} \sum_{i=1}^n \frac{1}{n} \sum_{j=1}^n \left\| \left[ \sum_{t=1}^T g_{i,t} (\mathbf{x}_{j,t}) \right]_+ \right\| = o(T).
\end{align}
\end{rem}

\section{Conclusion}
\label{sec_conc}
In this work, we considered a distributed OCO problem, with time-varying (potentially adversarial) constraints. We proposed a distributed primal-dual mirror descent based approach, in which the primal and dual updates are carried out locally at all the nodes. We utilized the challenging, but more realistic metric of dynamic regret and fit. Without assuming the more restrictive Slater's conditions, we achieved sublinear regret and fit under mild, commonly used assumptions. To the best of our knowledge, this is the first work to consider distributed OCO problem with non-coupled local time-varying constraints, and achieve sublinear dynamic regret and fit.
	
\bibliographystyle{IEEEtran}
\bibliography{abrv,References}

\newpage
\appendix
\section{Proof of Intermediate Lemmas}

\subsection{Proof of Lemma \ref{lemma_dual}}
\label{app_lemma_dual}
\begin{proof}[\unskip\nopunct]
We prove \eqref{eq_lemma_dual_1} by induction. $\mathbf{q}_{i,0} = \mathbf{0}$ by initialization. Also, since $g_{i,0}(\cdot) \equiv \mathbf{0}$, $\mathbf{b}_{i,1} = \mathbf{0}$. Therefore, $\mathbf{q}_{i,1} = \mathbf{0}_m \Rightarrow \| \mathbf{q}_{i,1} \| \leq \frac{F}{\beta_1}$. Let us assume \eqref{eq_lemma_dual_1} be true at time $t$, $\forall \ i \in [n]$. For $t+1$, first note that
\begin{align}
    (1 - \gamma_{t+1} \beta_{t+1}) \mathbf{q}_{i,t} + \gamma_{t+1} \mathbf{b}_{i,t+1} \leq (1 - \gamma_{t+1} \beta_{t+1}) \mathbf{q}_{i,t} + \gamma_{t+1} g_{i,t} (\mathbf{y}_{i,t+1}), \label{eq_alg_b_convex}
\end{align}
where \eqref{eq_alg_b_convex} follows from step \ref{alg_line_update_b} in Algorithm \ref{alg1} using convexity of $g_{i,t} (\cdot)$. Also, for vectors $\mathbf{x}, \mathbf{y}$ such that $\mathbf{x} \preceq \mathbf{y}$, $\| [\mathbf{x}]_+ \| \leq \| \mathbf{y} \|$. Therefore,
\begin{align}
    \left\| \mathbf{q}_{i,t+1} \right\| & \leq (1 - \gamma_{t+1} \beta_{t+1}) \left\| \mathbf{q}_{i,t} \right\| + \gamma_{t+1} \left\| g_{i,t} (\mathbf{y}_{i,t+1}) \right\| \nonumber \\
    & \leq (1 - \gamma_{t+1} \beta_{t+1}) \frac{F}{\beta_{t}} + \gamma_{t+1} F \nn \\
    & \overset{(a)}{\leq} \frac{F}{\beta_{t+1}}, \forall \ i \in [n],
\end{align}
where $(a)$ follows since $\{ \beta_t \}$ is a non-increasing sequence.

Next we prove \eqref{eq_lemma_dual_2}. First note that
\beq 
    \begin{aligned}
        \lnr \mathbf{q}_{i,t+1} - \mathbf{q}_{i} \rnr^2 & \overset{(b)}{\leq} \| ( 1 - \gamma_{t+1} \beta_{t+1}) \mathbf{q}_{i,t} + \gamma_{t+1} \mathbf{b}_{i,t+1} - \mathbf{q}_i \|^2 \\
        & = \| \mathbf{q}_{i,t} - \mathbf{q}_{i} \|^2 + \gamma_{t+1}^2 \| \mathbf{b}_{i,t+1} - \beta_{t+1} \mathbf{q}_{i,t} \|^2 \\
        & \quad + 2 \gamma_{t+1} \lp \mathbf{q}_{i,t} - \mathbf{q}_{i} \rp^T [\nabla g_{i,t} (\mathbf{x}_{i,t})] (\mathbf{y}_{i,t+1} - \mathbf{x}_{i,t}) \\
        & \quad + 2 \gamma_{t+1} (\mathbf{q}_{i,t} - \mathbf{q}_{i})^T g_{i,t} (\mathbf{x}_{i,t}) - 2 \gamma_{t+1} \beta_{t+1} (\mathbf{q}_{i,t} - \mathbf{q}_{i})^T \mathbf{q}_{i,t}
    \end{aligned}
    \label{eq_q_it_q_i}
\eeq
where $(b)$ follows from \eqref{eq_proj_ineq}. 
We first bound the second term in \eqref{eq_q_it_q_i}.
\begin{align}
    \| \mathbf{b}_{i,t+1} - \beta_{t+1} \mathbf{q}_{i,t} \| & \leq \| \mathbf{b}_{i,t+1} \| + \beta_{t+1} \| \mathbf{q}_{i,t} \| \nonumber \\
    & \overset{(c)}{\leq} \| \nabla g_{i,t}(\mathbf{x}_{i,t}) \| \| \mathbf{y}_{i,t+1} - \mathbf{x}_{i,t} \| + \| g_{i,t}(\mathbf{x}_{i,t}) \| + \beta_{t+1} \frac{F}{\beta_{t+1}} \nonumber \\
    & \overset{(d)}{\leq} G d(\mathcal(X)) + 2F = B_1, \label{eq_q_it_q_i_1}
\end{align}
where $(c)$ follows from Cauchy-Schwarz inequality and \eqref{eq_lemma_dual_1}. $(d)$ follows from assumptions (B3), (B4) and \eqref{eq_set_diameter}. For the fourth term in \eqref{eq_q_it_q_i},
\begin{align}
    - 2 \gamma_{t+1} \mathbf{q}_{i}^T [\nabla g_{i,t} (\mathbf{x}_{i,t})] (\mathbf{y}_{i,t+1} - \mathbf{x}_{i,t}) & \leq 2 \gamma_{t+1} \| \mathbf{q}_{i}\| \| \nabla g_{i,t} (\mathbf{x}_{i,t}) \| \|\mathbf{y}_{i,t+1} - \mathbf{x}_{i,t}\| \nonumber \\
    & \overset{(e)}{\leq} 2 \gamma_{t+1} \left( \frac{G^2 \alpha_{t+1}}{\mu} \| \mathbf{q}_{i} \|^2 + \frac{\mu}{4 \alpha_{t+1}} \| \mathbf{y}_{i,t+1} - \mathbf{x}_{i,t} \|^2 \right), \label{eq_q_it_q_i_2}
\end{align}
where $(e)$ follows from (B4) and using Young's inequality $xy \leq \frac{x^2}{2a} + \frac{a y^2}{2}$, for $a>0$. For the last term in \eqref{eq_q_it_q_i}
\begin{align}
    - 2 \gamma_{t+1} \beta_{t+1} (\mathbf{q}_{i,t} - \mathbf{q}_{i})^T \mathbf{q}_{i,t} & \overset{(f)}{=} - 2 \gamma_{t+1} \beta_{t+1} \lb \frac{\| \mathbf{q}_{i,t} - \mathbf{q}_{i} \|^2 + \| \mathbf{q}_{i,t} \|^2 - \| \mathbf{q}_{i} \|^2}{2} \rb \nonumber \\
    & \leq \gamma_{t+1} \beta_{t+1} \left( \| \mathbf{q}_{i} \|^2 - \| \mathbf{q}_{i,t} - \mathbf{q}_{i} \|^2 \right), \label{eq_q_it_q_i_3}
\end{align}
where $(f)$ follows from $\mathbf{x}^T \mathbf{y} = \frac{\| \mathbf{x} \|^2 + \| \mathbf{y} \|^2 - \| \mathbf{x} - \mathbf{y} \|^2}{2}$. Define 
$$\Delta_{t+1} \triangleq \sum_{i=1}^n \lb \| \mathbf{q}_{i,t+1} - \mathbf{q}_{i} \|^2 - (1 - \gamma_{t+1} \beta_{t+1}) \| \mathbf{q}_{i,t} - \mathbf{q}_{i} \|^2 \rb.$$
Using the upper bounds \eqref{eq_q_it_q_i_1}-\eqref{eq_q_it_q_i_3} in \eqref{eq_q_it_q_i}, and summing over all $i \in [n]$ we get \eqref{eq_lemma_dual_2}.
\end{proof}

% \newpage
\subsection{Proof of Lemma \ref{lemma_primal}}
\label{app_lemma_primal}
\begin{proof}[\unskip\nopunct]
Using convexity of $f_{i,t}$
\begin{align}
    f_{i,t} (\mathbf{x}_{i,t}) - f_{i,t} (\mathbf{x}^*_{t}) & \leq \left\langle \nabla f_{i,t} (\mathbf{x}_{i,t}), \mathbf{x}_{i,t} - \mathbf{x}^*_{t} \right\rangle \nonumber \\
    & = \left\langle \nabla f_{i,t} (\mathbf{x}_{i,t}), \mathbf{x}_{i,t} - \mathbf{y}_{i,t+1} \right\rangle + \left\langle \nabla f_{i,t} (\mathbf{x}_{i,t}), \mathbf{y}_{i,t+1} - \mathbf{x}^*_{t} \right\rangle \label{eq_lemma_dregret_one_term}
\end{align}
For the first term in \eqref{eq_lemma_dregret_one_term}
\begin{align}
    \left\langle \nabla f_{i,t} (\mathbf{x}_{i,t}), \mathbf{x}_{i,t} - \mathbf{y}_{i,t+1} \right\rangle & \leq G \| \mathbf{x}_{i,t} - \mathbf{y}_{i,t+1} \| \nonumber \\
    & \leq \frac{G^2 \alpha_{t+1}}{\mu} + \frac{\mu}{4 \alpha_{t+1}} \| \mathbf{y}_{i,t+1} - \mathbf{x}_{i,t} \|^2, \label{eq_lemma_dregret_one_term_1}
\end{align}
where \eqref{eq_lemma_dregret_one_term_1} follows from assumption (B4). For the second term in \eqref{eq_lemma_dregret_one_term}
\begin{align}
    \left\langle \nabla f_{i,t} (\mathbf{x}_{i,t}), \mathbf{y}_{i,t+1} - \mathbf{x}^*_{t} \right\rangle \leq \left\langle a_{i,t+1}, \mathbf{y}_{i,t+1} - \mathbf{x}^*_{t} \right\rangle + \left\langle \left[ \nabla g_{i,t}(\mathbf{x}_{i,t}) \right]^T \mathbf{q}_{i,t}, \mathbf{x}^*_{t} - \mathbf{y}_{i,t+1} \right\rangle, \label{eq_lemma_dregret_one_term_2}
\end{align}
where \eqref{eq_lemma_dregret_one_term_2} follows from step \ref{alg_line_update_a} in Algorithm \ref{alg1}. Next, we bound the two inner products in \eqref{eq_lemma_dregret_one_term_2}. First,
\beq
    \begin{aligned}
        & \left\langle \left[ \nabla g_{i,t}(\mathbf{x}_{i,t}) \right]^T \mathbf{q}_{i,t}, \mathbf{x}^*_{t} - \mathbf{x}_{i,t} + \mathbf{x}_{i,t} - \mathbf{y}_{i,t+1} \right\rangle \\
        & \qquad \qquad \leq \mathbf{q}_{i,t}^T \left[ \left( g_{i,t}(\mathbf{x}^*_{t}) - g_{i,t}(\mathbf{x}_{i,t}) \right) + \nabla g_{i,t}(\mathbf{x}_{i,t}) \left( \mathbf{x}_{i,t} - \mathbf{y}_{i,t+1} \right) \right],
    \end{aligned}
    \label{eq_lemma_dregret_one_term_2_1}
\eeq
where \eqref{eq_lemma_dregret_one_term_2_1} follows from the convexity of $g_{i,t} (\cdot)$. Secondly,
\begin{align}
    & \left\langle a_{i,t+1}, \mathbf{y}_{i,t+1} - \mathbf{x}^*_{t} \right\rangle \leq \frac{1}{\alpha_{t+1}} \Big[ \mathcal{D}_{\mathcal{R}} (\mathbf{x}^*_{t}, \mathbf{x}_{i,t}) - \mathcal{D}_{\mathcal{R}} (\mathbf{x}^*_{t}, \mathbf{y}_{i,t+1}) - \mathcal{D}_{\mathcal{R}} (\mathbf{y}_{i,t+1}, \mathbf{x}_{i,t}) \Big] \label{eq_lemma_dregret_one_term_2_2a} \\
    & \qquad \leq \frac{1}{\alpha_{t+1}} \mathcal{D}_{\mathcal{R}} (\mathbf{x}^*_{t}, \mathbf{x}_{i,t}) - \frac{1}{\alpha_{t+2}} \mathcal{D}_{\mathcal{R}} (\mathbf{x}^*_{t+1}, \mathbf{x}_{i,t+1}) + \frac{1}{\alpha_{t+2}} \mathcal{D}_{\mathcal{R}} (\mathbf{x}^*_{t+1}, \mathbf{x}_{i,t+1}) - \frac{1}{\alpha_{t+2}} \mathcal{D}_{\mathcal{R}} (\mathbf{x}^*_{t}, \mathbf{x}_{i,t+1}) \nonumber \\
    & \qquad \quad + \frac{1}{\alpha_{t+2}} \mathcal{D}_{\mathcal{R}} (\mathbf{x}^*_{t}, \mathbf{x}_{i,t+1}) - \frac{1}{\alpha_{t+2}} \mathcal{D}_{\mathcal{R}} (\mathbf{x}^*_{t}, \mathbf{y}_{i,t+1}) + \frac{1}{\alpha_{t+2}} \mathcal{D}_{\mathcal{R}} (\mathbf{x}^*_{t}, \mathbf{y}_{i,t+1}) - \frac{1}{\alpha_{t+1}} \mathcal{D}_{\mathcal{R}} (\mathbf{x}^*_{t}, \mathbf{y}_{i,t+1}) \nonumber \\
    & \qquad \quad - \frac{\mu}{2 \alpha_{t+1}} \| \mathbf{y}_{i,t+1} - \mathbf{x}_{i,t} \|^2, \label{eq_lemma_dregret_one_term_2_2b}
\end{align}
where, \eqref{eq_lemma_dregret_one_term_2_2a} follows from Lemma \ref{lemma_bregman}. \eqref{eq_lemma_dregret_one_term_2_2b} follows from \eqref{eq_bregman_lower_bd}. Next, we analyze some of the component terms in \eqref{eq_lemma_dregret_one_term_2_2b}.
\begin{align}
    & \frac{1}{\alpha_{t+2}} \left[ \mathcal{D}_{\mathcal{R}} (\mathbf{x}^*_{t+1}, \mathbf{x}_{i,t+1}) - \mathcal{D}_{\mathcal{R}} (\mathbf{x}^*_{t}, \mathbf{x}_{i,t+1}) \right] \leq \frac{K}{\alpha_{t+2}} \| \mathbf{x}^*_{t+1} - \mathbf{x}^*_{t} \|. \label{eq_lemma_dregret_one_term_2_2c}
\end{align}
which follows from assumption (C2). We use the upper bounds in \eqref{eq_lemma_dregret_one_term_2_1}, \eqref{eq_lemma_dregret_one_term_2_2b}-\eqref{eq_lemma_dregret_one_term_2_2c} to further upper bound \eqref{eq_lemma_dregret_one_term_2}. Next, we use \eqref{eq_lemma_dregret_one_term_1} and this new bound on \eqref{eq_lemma_dregret_one_term_2} to upper bound \eqref{eq_lemma_dregret_one_term}. Then, summing \eqref{eq_lemma_dregret_one_term} over $i \in [n]$ and $t \in [T]$, we get
\beq
    \begin{aligned}
        & \sum_{t=1}^T \sum_{i=1}^n \left[ f_{i,t} (\mathbf{x}_{i,t}) - f_{i,t} (\mathbf{x}^*_{t}) \right] \leq \frac{n G^2}{\mu} \sum_{t=1}^T \alpha_{t+1} + \sum_{t=1}^T \sum_{i=1}^n \frac{\mu}{4 \alpha_{t+1}} \| \mathbf{y}_{i,t+1} - \mathbf{x}_{i,t} \|^2 \\
        & \qquad + \sum_{t=1}^T \sum_{i=1}^n \mathbf{q}_{i,t}^T \left[ g_{i,t}(\mathbf{x}^*_{t}) - g_{i,t}(\mathbf{x}_{i,t}) \right]  + \sum_{t=1}^T \sum_{i=1}^n \mathbf{q}_{i,t}^T  \nabla g_{i,t}(\mathbf{x}_{i,t}) \left( \mathbf{x}_{i,t} - \mathbf{y}_{i,t+1} \right) \\
        & \qquad + \sum_{t=1}^T \sum_{i=1}^n \left[ \frac{1}{\alpha_{t+1}} \mathcal{D}_{\mathcal{R}} (\mathbf{x}^*_{t}, \mathbf{x}_{i,t}) - \frac{1}{\alpha_{t+2}} \mathcal{D}_{\mathcal{R}} (\mathbf{x}^*_{t+1}, \mathbf{x}_{i,t+1}) \right] \\
        & \qquad + \sum_{t=1}^T \frac{1}{\alpha_{t+2}} \left[\sum_{i=1}^n \mathcal{D}_{\mathcal{R}} (\mathbf{x}^*_{t}, \mathbf{x}_{i,t+1}) - \sum_{i=1}^n \mathcal{D}_{\mathcal{R}} (\mathbf{x}^*_{t}, \mathbf{y}_{i,t+1}) \right] + \sum_{t=1}^T \sum_{i=1}^n \frac{K}{\alpha_{t+2}} \| \mathbf{x}^*_{t+1} - \mathbf{x}^*_{t} \| \\
        & \qquad + \sum_{t=1}^T \sum_{i=1}^n \left( \frac{1}{\alpha_{t+2}} - \frac{1}{\alpha_{t+1}} \right) \mathcal{D}_{\mathcal{R}} (\mathbf{x}^*_{t}, \mathbf{y}_{i,t+1}) - \sum_{t=1}^T \sum_{i=1}^n \frac{\mu}{2 \alpha_{t+1}} \| \mathbf{y}_{i,t+1} - \mathbf{x}_{i,t} \|^2.
    \end{aligned}
    \label{eq_lemma_dregret_sum}
\eeq
Here, by the feasibility of optimal
\begin{equation}
    \mathbf{q}_{i,t}^T g_{i,t}(\mathbf{x}^*_{t}) \leq 0, \quad \forall \ i, t. \label{eq_q_times_g_zero}
\end{equation}
Also,
\begin{align}
    \sum_{i=1}^n \mathcal{D}_{\mathcal{R}} (\mathbf{x}^*_{t}, \mathbf{x}_{i,t+1}) - \sum_{i=1}^n \mathcal{D}_{\mathcal{R}} (\mathbf{x}^*_{t}, \mathbf{y}_{i,t+1}) & = \sum_{i=1}^n \Big[ \mathcal{D}_{\mathcal{R}} \big( \mathbf{x}^*_{t}, \sum_{j=1}^n [\mathbf{W}_t]_{ij} \mathbf{y}_{j,t+1} \big) - \mathcal{D}_{\mathcal{R}} (\mathbf{x}^*_{t}, \mathbf{y}_{i,t+1}) \Big] \nonumber \\
    & \leq \sum_{i=1}^n \Big[ \sum_{j=1}^n [\mathbf{W}_t]_{ij} \mathcal{D}_{\mathcal{R}} \big( \mathbf{x}^*_{t}, \mathbf{y}_{j,t+1} \big) - \mathcal{D}_{\mathcal{R}} (\mathbf{x}^*_{t}, \mathbf{y}_{i,t+1}) \Big] \label{eq_lemma_dregret_one_term_2_2d} \\
    & = \Big[ \sum_{j=1}^n \mathcal{D}_{\mathcal{R}} \big( \mathbf{x}^*_{t}, \mathbf{y}_{j,t+1} \big) - \sum_{i=1}^n \mathcal{D}_{\mathcal{R}} (\mathbf{x}^*_{t}, \mathbf{y}_{i,t+1}) \Big] = 0. \label{eq_lemma_dregret_one_term_2_2e}
\end{align}
\eqref{eq_lemma_dregret_one_term_2_2d} follows from step \ref{alg1_line_consensus} in Algorithm \ref{alg1}, and assumption (C1). \eqref{eq_lemma_dregret_one_term_2_2e} follows from assumption (A1) about the network. Finally,
\begin{align}
    \sum_{t=1}^T \sum_{i=1}^n \left( \frac{1}{\alpha_{t+2}} - \frac{1}{\alpha_{t+1}} \right) \mathcal{D}_{\mathcal{R}} (\mathbf{x}^*_{t}, \mathbf{y}_{i,t+1}) & \leq \sum_{t=1}^T \left( \frac{1}{\alpha_{t+2}} - \frac{1}{\alpha_{t+1}} \right) n K d(\mathcal(X)) \nn \\
    & \leq \frac{n K d(\mathcal(X))}{\alpha_{T+2}}, \label{eq_lemma_dregret_one_term_2_2f}
\end{align}
where \eqref{eq_lemma_dregret_one_term_2_2f} follows from \eqref{eq_bregman_bdd}, and using the fact that $\{ \alpha_t \}$ is a non-increasing sequence.

Using \eqref{eq_q_times_g_zero}, \eqref{eq_lemma_dregret_one_term_2_2e}-\eqref{eq_lemma_dregret_one_term_2_2f} in \eqref{eq_lemma_dregret_sum}, we get \eqref{eq_lemma_dregret}.
\end{proof}

% \newpage
\subsection{Proof of Lemma \ref{lemma_nw_error}}
\label{app_lemma_nw_error}
\begin{proof}[\unskip\nopunct]
We start with applying Lemma \ref{lemma_bregman} to step \ref{alg1_line_primal_update} in Algorithm \ref{alg1}.
\begin{align}
    \alpha_{t+1} \left\langle \mathbf{a}_{i,t+1}, \mathbf{y}_{i,t+1} - \mathbf{x}_{i,t} \right\rangle & \leq - \mathcal{D}_{\mathcal{R}} (\mathbf{x}_{i,t}, \mathbf{y}_{i,t+1}) - \mathcal{D}_{\mathcal{R}} (\mathbf{y}_{i,t+1}, \mathbf{x}_{i,t}) \nonumber \\
    \Rightarrow \mu \left\| \mathbf{y}_{i,t+1} - \mathbf{x}_{i,t} \right\|^2 & \leq \alpha_{t+1} \left\| \mathbf{a}_{i,t+1} \right\| \left\| \mathbf{y}_{i,t+1} - \mathbf{x}_{i,t} \right\| \label{eq_nw_error_1}
\end{align}
where \eqref{eq_nw_error_1} follows from \eqref{eq_bregman_lower_bd}. Consequently,
\begin{align}
    \left\| \mathbf{y}_{i,t+1} - \mathbf{x}_{i,t} \right\| & \leq \frac{\alpha_{t+1}}{\mu} \left\| \mathbf{a}_{i,t+1} \right\| \leq \frac{\alpha_{t+1}}{\mu} \left( \left\| \nabla f_{i,t}(\mathbf{x}_{i,t}) \right\| + \left\| \nabla g_{i,t}(\mathbf{x}_{i,t}) \mathbf{q}_{i,t} \right\| \right) \nonumber \\
    & \leq \frac{\alpha_{t+1}}{\mu} \left( G + G \frac{F}{\beta_t} \right) \leq \frac{G \alpha_{t+1}}{\mu} \left( 1 + \frac{F}{\beta_{t+1}} \right). \label{eq_nw_error_2}
\end{align}
Note that the bound in \eqref{eq_nw_error_2} is independent of $i$. Next, define $\mathbf{e}_{i,t} = \mathbf{y}_{i,t+1} - \mathbf{x}_{i,t}$. Note that,
\begin{align}
    \mathbf{y}_{i,t+1} = \mathbf{x}_{i,t} + \mathbf{e}_{i,t} = \sum_{j=1}^n [\mathbf{W}_t]_{ij} \mathbf{y}_{j,t} + \mathbf{e}_{i,t}, \label{eq_nw_error_y_i}
\end{align}
and
\begin{align}
    \bar{\mathbf{y}}_{t+1} &= \bar{\mathbf{x}}_{t} + \bar{\mathbf{e}}_{t} = \frac{1}{n} \sum_{i=1}^n \sum_{j=1}^n [\mathbf{W}_t]_{ij} \mathbf{y}_{j,t} + \bar{\mathbf{e}}_{t} \nonumber \\
    & = \bar{\mathbf{y}}_{t} + \bar{\mathbf{e}}_{t} = \sum_{\tau=0}^t \bar{\mathbf{e}}_{\tau} \label{eq_nw_error_y_avg}.
\end{align}
Since we assume $f_{i,0}  (\cdot) \equiv 0, \forall \ i \in [n]$, $g_{i,0} (\cdot) = \mathbf{0}, \forall \ i \in [n]$ and $\mathbf{x}_{i,0} = \mathbf{0}$, from the algorithm, $\mathbf{y}_{i,1} = \mathbf{0}, \forall \ i \in [n]$. Further, defining
\begin{align}
    \mathbf{y}_t \triangleq \left[ \mathbf{y}_{1,t}^T, \hdots, \mathbf{y}_{n,t}^T \right]^T, \ \mathbf{e}_t \triangleq \left[ \mathbf{e}_{1,t}^T, \hdots, \mathbf{e}_{n,t}^T \right]^T
\end{align}
we see from \eqref{eq_nw_error_y_i} that $\mathbf{y}_{t+1} = (\mathbf{W} \otimes I) \mathbf{y}_t + (I \otimes I)  \mathbf{e}_t$. Hence
\begin{align}
    \mathbf{y}_{i,t+1} = \sum_{\tau = 0}^t \sum_{j=1}^n \left[ \mathbf{W}^{t-\tau} \right]_{ij} \mathbf{e}_{j, \tau}. \label{eq_nw_error_y_i_wrt_e}
\end{align}
Substituting \eqref{eq_nw_error_y_i_wrt_e} in \eqref{eq_nw_error_y_i}, and using it with \eqref{eq_nw_error_y_avg}, we
\begin{align}
    \mathbf{x}_{i,t} - \bar{\mathbf{x}}_t &= \sum_{j=1}^n [\mathbf{W}]_{ij} \mathbf{y}_{j,t} - \bar{\mathbf{y}}_t \nonumber \\
    &= \sum_{j=1}^n \left( [\mathbf{W}]_{ij} - \frac{1}{n} \right) \sum_{\tau = 0}^{t-1} \sum_{k=1}^n \left[ \mathbf{W}^{t-1-\tau} \right]_{jk} \mathbf{e}_{k, \tau} \nonumber \\
    &= \sum_{\tau = 0}^{t-1} \sum_{k=1}^n \left( \left[ \mathbf{W}^{t-\tau} \right]_{ik} - \frac{1}{n} \right) \mathbf{e}_{k, \tau}
\end{align}
Consequently,
\begin{align}
    \left\| \mathbf{x}_{i,t} - \bar{\mathbf{x}}_t \right\| & \leq \sum_{\tau = 0}^{t-1} \sum_{k=1}^n \left| \left[ \mathbf{W}^{t-\tau} \right]_{ik} - \frac{1}{n} \right| \left\| \mathbf{e}_{k, \tau} \right\|, \label{eq_nw_error_3}
\end{align}
A standard property of doubly-stochastic matrices is that
\begin{equation}
    \label{eq_doubly_stoch_err}
    \left| \left[ \mathbf{W}^{t} \right]_{ik} - \frac{1}{n} \right| \leq \sqrt{n} \sigma_2^t(\mathbf{W}).
\end{equation}
Substituting the bound on $\| \mathbf{e}_{k, \tau} \|$ from \eqref{eq_nw_error_2}, and \eqref{eq_doubly_stoch_err} we get \eqref{eq_nw_error}.
\end{proof}

% \newpage
\subsection{Proof of Lemma \ref{lemma_regret_fit}}
\label{app_lemma_regret_fit}
\begin{proof}[\unskip\nopunct]
Substituting $\mathbf{q}_i = \mathbf{0}$, $\forall \ i \in [n]$ in \eqref{eq_bd_gc_q}, followed by adding \eqref{eq_bd_gc_q} and \eqref{eq_lemma_dregret}, we get \eqref{eq_lemma_dregret_2}. Similarly, substituting
\begin{align}
    \bar{\mathbf{q}}_i = \frac{\left[ \sum_{t=1}^T g_{i,t} (\mathbf{x}_{i,t}) \right]_+}{2\left[ \frac{1}{2 \gamma_1} + \sum_{t=1}^T \left( \frac{G^2 \alpha_{t+1}}{\mu} + \frac{\beta_{t+1}}{2} \right) \right]} \label{eq_q_for_fit}
\end{align}
% in \eqref{eq_bd_gc_q}, followed by adding \eqref{eq_bd_gc_q} and \eqref{eq_lemma_dregret}, we get \eqref{eq_lemma_dfit}.
% \end{proof}

% \begin{proof}[\unskip\nopunct]
% If instead, we substitute $\bar{\mathbf{q}}_i$ from \eqref{eq_q_for_fit}
and $\bar{\mathbf{q}}_j = \mathbf{0}$, $\forall \ j \neq i$, in \eqref{eq_bd_gc_q}. Subsequently, by adding \eqref{eq_bd_gc_q} and \eqref{eq_lemma_dregret}, we get
\beq
    \begin{aligned}
        & \sum_{i=1}^n \left\| \left[ \sum_{t=1}^T g_{i,t} (\mathbf{x}_{i,t}) \right]_+ \right\|^2 \leq 4\left[ \frac{1}{2 \gamma_1} + \sum_{t=1}^T \left( \frac{G^2 \alpha_{t+1}}{\mu} + \frac{\beta_{t+1}}{2} \right) \right] \Bigg\{ 2nFT \\
        & \qquad + \frac{n B_1^2}{2} \sum_{t=1}^T \gamma_{t+1} + \frac{n G^2}{\mu} \sum_{t=1}^T \alpha_{t+1} + \sum_{i=1}^n \left[ \frac{1}{\alpha_{2}} \mathcal{D}_{\mathcal{R}} (\mathbf{x}^*_{1}, \mathbf{x}_{i,1}) - \frac{1}{\alpha_{T+2}} \mathcal{D}_{\mathcal{R}} (\mathbf{x}^*_{T+1}, \mathbf{x}_{i,T+1}) \right] \\
        & \qquad + \frac{n K}{\alpha_{T+2}} \sum_{t=1}^T \| \mathbf{x}^*_{t+1} - \mathbf{x}^*_{t} \| + \frac{n K d(\mathcal(X))}{\alpha_{T+2}} - \frac{1}{2} \sum_{t=1}^T \left( \frac{1}{\gamma_t} - \frac{1}{\gamma_{t+1}} + \beta_{t+1} \right) \sum_{i=1}^n \| \mathbf{q}_{i,t} - \bar{\mathbf{q}}_i \|^2 \Bigg\}.
    \end{aligned}
    \label{eq_lemma_dfit_one_node}
\eeq
\end{proof}

\newpage
\section{Proof of Theorem \ref{thm_regret_fit}}
\label{app_thm_regret_fit}
\begin{proof}[\unskip\nopunct]
We begin with one of the constituent terms in \eqref{eq_defn_dregret}.
\begin{align}
    f_t(\mathbf{x}_{i,t}) - f_t(\mathbf{x}^*_t) &= f_t(\mathbf{x}_{i,t}) - f_t(\bar{\mathbf{x}}_{t}) + f_t(\bar{\mathbf{x}}_{t}) - f_t(\mathbf{x}^*_t) \nonumber \\
    & \leq L \left\| \mathbf{x}_{i,t} - \bar{\mathbf{x}}_{t} \right\| + f_t(\bar{\mathbf{x}}_{t}) - f_t(\mathbf{x}^*_t) \\
    & = \frac{1}{n} \sum_{j=1}^n \left\{ f_{j,t} (\bar{\mathbf{x}}_{t}) - f_{j,t} (\mathbf{x}^*_t) + f_{j,t} (\mathbf{x}_{j,t}) - f_{j,t} (\mathbf{x}_{j,t}) \right\} + L \left\| \mathbf{x}_{i,t} - \bar{\mathbf{x}}_{t} \right\| \nonumber \\
    & \leq \frac{1}{n} \sum_{j=1}^n \left\{ f_{j,t} (\mathbf{x}_{j,t}) - f_{j,t} (\mathbf{x}^*_t) \right\} + L \left\| \mathbf{x}_{i,t} - \bar{\mathbf{x}}_{t} \right\| + \frac{L}{n} \sum_{j=1}^n \left\| \mathbf{x}_{j,t} - \bar{\mathbf{x}}_{t} \right\|.
\end{align}
We use assumption (B2) to get both \eqref{eq_use_f_lipschitz_1}, \eqref{eq_use_f_lipschitz_2}.
Now, from the definition of dynamic regret \eqref{eq_defn_dregret}, we get
\begin{align}
    \mathbf{Reg}_T^d & \leq \frac{1}{n} \sum_{i=1}^n \sum_{t=1}^T \frac{1}{n} \sum_{j=1}^n \left\{ f_{j,t} (\mathbf{x}_{j,t}) - f_{j,t} (\mathbf{x}^*_t) \right\} + \frac{2L}{n} \sum_{i=1}^n \sum_{t=1}^T \left\| \mathbf{x}_{i,t} - \bar{\mathbf{x}}_{t} \right\|. \label{eq_bd_dregret_1}
\end{align}
We use Lemma \ref{lemma_regret_fit} to bound the first term, and Lemma \ref{lemma_nw_error} to bound the second term in \eqref{eq_bd_dregret_1}. But first, we bound all the terms in the upper bounds in Lemma \ref{lemma_nw_error} and \ref{lemma_regret_fit}, using the step-sizes in \eqref{eq_step_size_choice}.

For a constant $p < 1$ and $T \in \mathbb{N}_+$
\begin{align}
    \sum_{t=1}^T \frac{1}{t^p} \leq 1 + \int_1^T \frac{1}{t^p} dt = \frac{T^{1-p} - 1}{1-p} + 1 \leq \frac{T^{1-p}}{1-p}.
\end{align}
Hence,
\begin{align}
    & \sum_{t=1}^T \gamma_{t+1} \leq \frac{T^{b}}{b}, \ \sum_{t=1}^T \alpha_{t+1} \leq \frac{T^{1-a}}{1-a}, \ \sum_{t+1}^T \beta_{t+1} \leq \frac{T^{1-b}}{1-b}, \label{eq_sum_step_size} \\
    & \frac{1}{\alpha_{2}} \mathcal{D}_{\mathcal{R}} (\mathbf{x}^*_{1}, \mathbf{x}_{i,1}) - \frac{1}{\alpha_{T+2}} \mathcal{D}_{\mathcal{R}} (\mathbf{x}^*_{T+1}, \mathbf{x}_{i,T+1}) \leq 2 K d(\mathcal(X)), \label{eq_Dr_Dr_bound} \\
    & \frac{1}{\alpha_{T+2}} = (T+2)^a \leq 2T^a, \qquad \forall \ T \geq 2. \label{eq_final_step_size}
\end{align}
\eqref{eq_Dr_Dr_bound} follows from $\mathcal{D}_{\mathcal{R}} (\cdot, \cdot) \geq 0$, \eqref{eq_bregman_bdd} and $\frac{1}{\alpha_2} = 2^a < 2$. Also, $\left( \frac{1}{\gamma_t} - \frac{1}{\gamma_{t+1}} + \beta_{t+1} \right) > 0$ follows from \eqref{eq_step_size_choice}. Therefore, we ignore the last term in both \eqref{eq_lemma_dregret_2} and \eqref{eq_lemma_dfit}. Also, substituting \eqref{eq_nw_error} in the second term in \eqref{eq_bd_dregret_1}, we get
\begin{align}
    2L \sum_{t=1}^T \sum_{\tau = 0}^{t-1} \sqrt{n} \sigma_2^{t-\tau} (\mathbf{W}) \frac{G \alpha_{\tau+1}}{\mu} \left( 1 + \frac{F}{\beta_{\tau+1}} \right) &= \frac{2 L G \sqrt{n}}{\mu} \sum_{t=1}^T \sum_{\tau = 0}^{t-1} \left( \frac{\sigma_2^{t-\tau}}{(\tau+1)^a} + F \frac{\sigma_2^{t-\tau}}{(\tau+1)^{a-b}} \right) \nonumber \\
    &= \frac{2 L G \sqrt{n}}{\mu} \sum_{t=1}^T \sum_{\tau = 1}^{t} \left( \frac{\sigma_2^{t+1-\tau}}{\tau^a} + F \frac{\sigma_2^{t+1-\tau}}{\tau^{a-b}} \right) \nonumber \\
    & \leq \frac{2 L G \sqrt{n} \sigma_2(\mathbf{W})}{\mu (1 - \sigma_2(\mathbf{W}))} \sum_{t=1}^T \left( \frac{1}{\tau^a} + F \frac{1}{\tau^{a-b}} \right) \label{eq_total_nw_error_inter} \\
    & = \frac{2 L G \sqrt{n} \sigma_2(\mathbf{W})}{\mu (1 - \sigma_2(\mathbf{W}))} \left( \frac{T^{1-a}}{1-a} + \frac{F T^{1+b-a}}{1-a+b} \right). \label{eq_total_nw_error}
\end{align}
where in \eqref{eq_total_nw_error_inter}, we use the following inequality. Given a non-negative sequence $\{\delta_t \}$, and $0 < \lambda < 1$, it holds
\begin{equation}
    \sum_{t=1}^T \sum_{\tau = 1}^t \delta_{\tau + 1} \lambda^{t-\tau} = \sum_{t=1}^T \delta_{t + 1} \sum_{\tau = 0}^{T - \tau} \lambda^{\tau} \leq \frac{1}{1-\lambda} \sum_{t=1}^T \delta_{t + 1}.
\end{equation}
Using \eqref{eq_sum_step_size}-\eqref{eq_final_step_size}, \eqref{eq_total_nw_error}, and Lemma \ref{lemma_nw_error}, \ref{lemma_regret_fit} we upper bound the dynamic regret \eqref{eq_bd_dregret_1} as
\begin{align}
    \mathbf{Reg}_T^d & \leq \frac{B_1^2}{2} \frac{T^b}{b} + \frac{G^2}{\mu} \frac{T^{1-a}}{1-a} + 2 T^a K d(\mathcal(X)) + 2 K T^a C_T^* \nonumber \\
    & \quad + 2 K d(\mathcal(X)) + \frac{2 L G \sqrt{n} \sigma_2(\mathbf{W})}{\mu (1 - \sigma_2(\mathbf{W}))} \left( \frac{T^{1-a}}{1-a} + \frac{F T^{1+b-a}}{1-a+b} \right) \label{eq_bd_dregret_2}
\end{align}
Since $a > b$, for large enough $T$, $T^a$ dominates $T^b$, while, $T^{1-a+b}$ will dominate $T^{1-a}$.

To upper bound the fit, we use \eqref{eq_sum_step_size}-\eqref{eq_final_step_size}, and Lemma \ref{lemma_regret_fit}.
\begin{align}
    \frac{1}{n} \sum_{i=1}^n \left\| \left[ \sum_{t=1}^T g_{i,t} (\mathbf{x}_{i,t}) \right]_+ \right\|^2 & \leq 4 \left[ \frac{1}{2 \gamma_1} + \frac{G^2}{\mu} \frac{T^{1-a}}{1-a} + \frac{T^{1-b}}{2(1-b)} \right] \Bigg\{ 2FT + 2 K d(\mathcal(X)) \nonumber \\
    & \quad + \frac{B_1^2}{2} \frac{T^b}{b} + \frac{G^2}{\mu} \frac{T^{1-a}}{1-a} + 2 T^a K d(\mathcal(X)) + 2 K T^a C_T^* \Bigg\}.
\end{align}
For large enough $T$, $T^{2-b}$ dominates $T, T^{2-a}$, $T^{2-a-b}$ and $T^{1+a-b}$. Also, $T^{1+a-b}$ dominates $T$.
\end{proof}

\subsection{Fixed Step Sizes}
Assume $\alpha_t = \frac{1}{T^a}, \beta_t = \frac{1}{T^b}, \gamma_t = \frac{1}{T^c}$. Here, $a>b$. Then, the individual terms in Lemma \ref{lemma_regret_fit} can be bounded as follows.
\beq
    \begin{aligned}
        & \sum_{t=1}^T \gamma_{t+1} = T^{1-c}, \qquad \sum_{t=1}^T \alpha_{t+1} = T^{1-a}, \qquad \sum_{t+1}^T \beta_{t+1} = T^{1-b}, \\
        & \left[ \frac{1}{\alpha_{2}} \mathcal{D}_{\mathcal{R}} (\mathbf{x}^*_{1}, \mathbf{x}_{i,1}) - \frac{1}{\alpha_{T+2}} \mathcal{D}_{\mathcal{R}} (\mathbf{x}^*_{T+1}, \mathbf{x}_{i,T+1}) \right] \leq T^a K d(\mathcal(X)), \\
        & \sum_{t=1}^T \frac{n K}{\alpha_{t+2}} \| \mathbf{x}^*_{t+1} - \mathbf{x}^*_{t} \| \leq n K T^a \sum_{t=1}^T \| \mathbf{x}^*_{t+1} - \mathbf{x}^*_{t} \|, \\
        & \frac{n K d(\mathcal(X))}{\alpha_{T+2}} \leq n K d(\mathcal(X)) T^a.
    \end{aligned}
\eeq
Similarly, we can bound the individual terms in Lemma \ref{lemma_nw_error}.
\begin{align}
    2L \sum_{t=1}^T \sum_{\tau = 0}^{t-1} \sqrt{n} \sigma_2^{t-\tau} (\mathbf{W}) \frac{G \alpha_{\tau+1}}{\mu} \left( 1 + \frac{F}{\beta_{\tau+1}} \right) &\leq \frac{2 L G \sqrt{n}}{\mu} \left( T^{-a} + F T^{b-a} \right) \frac{T}{1 - \sigma_2(\mathbf{W})} \nonumber \\
    &= \frac{2 L G \sqrt{n}}{\mu (1 - \sigma_2(\mathbf{W}))} \left( T^{1-a} + F T^{1+b-a} \right).
\end{align}
Henceforth, we use denote the cumulative consecutive difference of optimal points $C_T^* \triangleq \sum_{t=1}^T \| \mathbf{x}^*_{t+1} - \mathbf{x}^*_{t} \|$.
Consequently, we can upper bound the dynamic regret as
\begin{align}
    \mathbf{Reg}_T^d & \leq \frac{1}{n} \sum_{i=1}^n \sum_{t=1}^T \left\{ f_{i,t} (\mathbf{x}_{i,t}) - f_{i,t} (\mathbf{x}^*_t) \right\} + \frac{2L}{n} \sum_{i=1}^n \sum_{t=1}^T \left\| \mathbf{x}_{i,t} - \bar{\mathbf{x}}_{t} \right\| \nonumber \\
    & \leq \frac{B_1^2}{2} T^{1-c} + \frac{G^2}{\mu} T^{1-a} + 2 T^a K d(\mathcal(X)) + K T^a C_T^* + \frac{2 L G \sqrt{n}}{\mu (1 - \sigma_2(\mathbf{W}))} \left( T^{1-a} + F T^{1+b-a} \right) \\
    & \leq A_{1} T^{\max \{ a, 1-c, 1-a+b \}} + K T^a C_T^*
\end{align}

And,
\begin{align}
    \frac{1}{n} \sum_{i=1}^n \left\| \left[ \sum_{t=1}^T g_{i,t} (\mathbf{x}_{i,t}) \right]_+ \right\|^2 & \leq 4\left[ \frac{1}{2 \gamma_1} + \frac{G^2}{\mu} T^{1-a} + \frac{T^{1-b}}{2} \right] \big\{ 2FT \nonumber \\
    & \quad + \frac{B_1^2}{2} T^{1-c} + \frac{G^2}{\mu} T^{1-a}  + 2 T^a K d(\mathcal(X)) + K T^a C_T^* \big\} \\
    & \leq C_1 T^{\max \{2-a, 2-b\}} + C_2 T^{\max\{ 1, 1+a-b \}} C_T^* \\
    & = C_1 T^{2-b} + C_2 T^{\max\{ 1, 1+a-b \}} C_T^*
\end{align}

\end{document}